\documentclass[11pt]{article}
\usepackage[usenames,dvipsnames,svgnames]{xcolor}
\usepackage[shortlabels]{enumitem}
\usepackage[framemethod=TikZ]{mdframed}
\usepackage{amsmath,amssymb,amsthm}
\usepackage{epigraph}
\usepackage{caption}
\usepackage[colorlinks]{hyperref}
\usepackage{microtype}
\usepackage{mathtools}
\usepackage[headsepline]{scrlayer-scrpage}
\usepackage{thmtools}
\usepackage{listings}
\usepackage{float}
\usepackage[section]{placeins}
\usepackage{derivative}
\usepackage[margin=1in]{geometry}

\providecommand{\clubg}[1]{\bgroup\color{green!40!black}[$#1\clubsuit$]\egroup}

\DeclareMathOperator{\Av}{Av}
\DeclareMathOperator{\rev}{rev}

\reversemarginpar

\mdfdefinestyle{mdbluebox}{roundcorner=10pt,innerbottommargin=9pt,
    linecolor=blue,backgroundcolor=TealBlue!5,}
\declaretheoremstyle[headfont=\sffamily\bfseries\color{MidnightBlue},
    mdframed={style=mdbluebox},]{thmbluebox}
\mdfdefinestyle{mdredbox}{frametitlefont=\bfseries,innerbottommargin=8pt,
    nobreak=true,backgroundcolor=Salmon!5,linecolor=RawSienna,}
\declaretheoremstyle[headfont=\bfseries\color{RawSienna},
    mdframed={style=mdredbox},headpunct={\\[3pt]},postheadspace=0pt,]{thmredbox}
\mdfdefinestyle{mdgreenbox}{linecolor=ForestGreen,backgroundcolor=ForestGreen!5,
    linewidth=2pt,rightline=false,leftline=true,topline=false,bottomline=false,}
\declaretheoremstyle[headfont=\bfseries\sffamily\color{ForestGreen!70!black},
    mdframed={style=mdgreenbox},headpunct={ --- },]{thmgreenbox}
\mdfdefinestyle{mdblackbox}{linecolor=black,backgroundcolor=RedViolet!5!gray!5,
    linewidth=3pt,rightline=false,leftline=true,topline=false,bottomline=false,}
\declaretheoremstyle[mdframed={style=mdblackbox}]{thmblackbox}

\declaretheorem[style=thmbluebox,name=Theorem,numbered=no]{theorem*}
\declaretheorem[style=thmbluebox,name=Lemma,numbered=no]{lemma*}

\declaretheorem[style=thmgreenbox,name=Claim,numbered=no]{claim*}

\declaretheorem[style=thmblackbox,name=Remark,numbered=no]{remark*}

\declaretheorem[style=thmgreenbox,name=Definition,numbered=no]{definition*}

\declaretheorem[style=thmblackbox,name=Example,numbered=no]{example*}

\newlist{walk}{enumerate}{3}
\setlist[walk]{label=\bfseries (\alph*)}
\newenvironment{subproof}[1][Proof]{%
\begin{proof}[#1] }%
{\end{proof}}

\usepackage{asymptote}
\begin{asydef}
size(10cm); // set a reasonable default
usepackage("amsmath");
usepackage("amssymb");
settings.tex="pdflatex";
settings.outformat="pdf";
import geometry;
void filldraw(picture pic = currentpicture, conic g, pen fillpen=defaultpen, pen drawpen=defaultpen) { filldraw(pic, (path) g, fillpen, drawpen); }
void fill(picture pic = currentpicture, conic g, pen p=defaultpen) { filldraw(pic, (path) g, p); }
pair foot(pair P, pair A, pair B) { return foot(triangle(A,B,P).VC); }
pair centroid(pair A, pair B, pair C) { return (A+B+C)/3; }
\end{asydef}

\begin{document}
\title{Computational Approach to the $SC_{231}$ Consecutive-Pattern-Avoiding Stack Sort}
\author{Kai Yi}
\date{\today}
\maketitle

\begin{abstract}
Defant and Zheng introduced a consecutive-pattern-avoiding stack sort map $SC_{\sigma}$, where the stack must avoid a consecutive pattern $\sigma$. Seidel and Sun disproved a conjecture in Defant and Zheng's paper about the maximum sort-number of a length $n$ permutation under $SC_{231}$. In this paper, we compute sort-numbers for each permutation of length up to $14$, and we estimate the average sort-numbers up to length $1000$. Our results suggest the maximum and average sort-numbers grow faster than linear with respect to $n$ for the tested ranges, though the long-term behavior remains unclear. We also prove properties of $SC_{231}$ mathematically, such as a $n-1$ lower bound and a $\frac{(n+1)(n-2)}{2}$ upper bound for the maximum sort-number of length $n$ permutations.
\end{abstract}

\section{Introduction and Preliminaries}

Stack sorting was first introduced in 1968 by Knuth in \cite{Knuth}, who proved a permutation can be sorted into the identity permutation using the stack sort machine (push and pop operations) if and only if it avoids the pattern $231$. The formal definition of pattern avoidance is:

Let $p$ and $q$ be permutations. Let the entries of $q$ be $q_1q_2 \dots q_k$. If there exists a subsequence of entries $p_1 p_2 \dots p_k$ in $p$ such that $p_i < p_j$ if and only if $q_i < q_j$ for all indices $i, j$, we say $p$ contains the pattern $q$. Otherwise, we say $p$ avoids the pattern $q$.

For example, $13425$ avoids the pattern $321$ but contains the pattern $123$ because of the entries $1, 3, 5$.
\newline

In 1990, West \cite{West} introduced the stack sort function $s : S_n \rightarrow S_n$, which "pushes" an entry into a stack if the entries in the stack would remain in increasing order from top to bottom and "pops" an entry from the top of the stack if this would not happen. This is similar to the Tower of Hanoi, where you can not place a large disk on top of a small one.

For example, here is how the permutation $15243$ is sorted using West's function:
\begin{figure}[H]
\centering
\begin{asy}
draw((0,0)--(2,0)--(2,-2)--(3,-2)--(3,0)--(5,0));
draw((5,-1)--(6,-1), arrow=Arrow());
draw((6,0)--(8,0)--(8,-2)--(9,-2)--(9,0)--(11,0));
draw((11,-1)--(12,-1), arrow=Arrow());
draw((12,0)--(14,0)--(14,-2)--(15,-2)--(15,0)--(17,0));
draw((17,-1)--(18,-1), arrow=Arrow());
draw((18,0)--(20,0)--(20,-2)--(21,-2)--(21,0)--(23,0));

draw((-1,-4)--(0,-4), arrow=Arrow());
draw((0,-3)--(2,-3)--(2,-5)--(3,-5)--(3,-3)--(5,-3));
draw((5,-4)--(6,-4), arrow=Arrow());
draw((6,-3)--(8,-3)--(8,-5)--(9,-5)--(9,-3)--(11,-3));
draw((11,-4)--(12,-4), arrow=Arrow());
draw((12,-3)--(14,-3)--(14,-5)--(15,-5)--(15,-3)--(17,-3));
draw((17,-4)--(18,-4), arrow=Arrow());
draw((18,-3)--(20,-3)--(20,-5)--(21,-5)--(21,-3)--(23,-3));

draw((-1,-7)--(0,-7), arrow=Arrow());
draw((0,-6)--(2,-6)--(2,-8)--(3,-8)--(3,-6)--(5,-6));
draw((5,-7)--(6,-7), arrow=Arrow());
draw((6,-6)--(8,-6)--(8,-8)--(9,-8)--(9,-6)--(11,-6));
draw((11,-7)--(12,-7), arrow=Arrow());
draw((12,-6)--(14,-6)--(14,-8)--(15,-8)--(15,-6)--(17,-6));

label("\footnotesize15243", (4, 0.5));
label("\footnotesize5243", (10, 0.5));
label("\footnotesize5243", (16, 0.5));
label("\footnotesize243", (22, 0.5));
label("\footnotesize43", (4, -2.5));
label("\footnotesize43", (10, -2.5));
label("\footnotesize3", (16, -2.5));

label("\footnotesize1", (13, 0.5));
label("\footnotesize1", (19, 0.5));
label("\footnotesize1", (1, -2.5));
label("\footnotesize12", (7, -2.5));
label("\footnotesize12", (13, -2.5));
label("\footnotesize12", (19, -2.5));
label("\footnotesize123", (1, -5.5));
label("\footnotesize1234", (7, -5.5));
label("\footnotesize12345", (13, -5.5));

label("\footnotesize1", (8.5, -1.5));
label("\footnotesize5", (20.5, -1.5));
label("\footnotesize5", (2.5, -4.5));
label("\footnotesize2", (2.5, -3.8));
label("\footnotesize5", (8.5, -4.5));
label("\footnotesize5", (14.5, -4.5));
label("\footnotesize4", (14.5, -3.8));
label("\footnotesize5", (20.5, -4.5));
label("\footnotesize4", (20.5, -3.8));
label("\footnotesize3", (20.5, -3.1));
label("\footnotesize5", (2.5, -7.5));
label("\footnotesize4", (2.5, -6.8));
label("\footnotesize5", (8.5, -7.5));
\end{asy}
\caption*{\footnotesize Figure 1: Stack Sorting $15243$ using West's function}
\end{figure}

In the second step, we pop $1$ rather than push $5$ because pushing $5$ would result in $5$ (a large entry) being in the stack on top of $1$ (a small entry).
\newline

Since West's 1990 dissertation, many have explored variants of the original stack sort function. For example, in 2020, Cerbai, Claesson, and Ferrari \cite{CCF} introduced the pattern-avoiding stack sorting function, $s_{\sigma}$. This stack differs from West's stack because the numbers in the stack, arranged top to bottom, must avoid the pattern $\sigma$. In this new definition, West's stack sorting function is $s_{21}$.

A few months later, Defant and Zheng \cite{DefantZheng} introduce a related variant, the stack sorting function $SC_{\sigma}$. Unlike $s_{\sigma}$, the stack must avoid a consecutive pattern instead of a regular pattern.

For example, $15243$ contains the pattern $123$ because of entries $1, 2, 3$, but it avoids the consecutive pattern $\underline{123}$ because there does not exist three consecutive entries in increasing order. In our paper, we will underline a consecutive pattern to denote that it is consecutive.
\newline

We will look at $SC_{231}$ specifically. Let $\Av_n(\sigma_1, \sigma_2, \dots, \sigma_k)$ denote the set of permutations in $S_n$ that avoid the patterns $\sigma_1, \sigma_2, \dots, \sigma_k$. In Defant and Zheng's paper \cite{DefantZheng}, the authors proposed the conjecture:

\textbf{Conjecture 1.1:} Let $n \geq 3$. We have $SC_{231}^{2n-4}(\pi) \in \Av_n(132, 231)$ for every $\pi \in S_n$. Furthermore, there exists a permutation $\tau \in S_n$ such that $SC_{231}^{2n-5}(\tau) \not\in \Av_n(132, 231)$.

The motivation behind this conjecture is the fact that the periodic points of $SC_{\sigma}$ are the permutations in $\Av_n(\underline{\sigma}, \underline{\rev(\sigma)})$. This was proved in Defant and Zheng's paper \cite{DefantZheng} for patterns $\sigma \in S_3$, then proved in Siedel and Sun's paper \cite{SiedelSun} for $\sigma \in S_k$, $k \geq 3$. Note that $\Av_n(\underline{231}, \underline{\rev(231)})$ is the same as $\Av_n(231, \rev(231))$ because both sets represent length $n$ permutations that have no peaks.

This conjecture was first verified in \cite{DefantZheng} up to $n = 9$, then disproved in \cite{SiedelSun} using the counterexample $(4, 6, 8, 5, 11, 7, 2, 9, 10, 3, 1) \in S_{11}$. This leads to the question: how many times do we really need to apply $SC_{231}$ to force a length $11$ permutation into a periodic point? We will explore this in more depth using a computational approach.

Another paper that uses $SC_{231}$ is \cite{Zhao}, in which they denote $SC_{231}$ as $SC_{2\underline{31}}$. The paper \cite{Kemeklis} also explores $SC_{231}$, finding the fertility numbers of specific examples to prove that for any positive number $f$, there exists a permutation $p$ such that $|SC_{231}^{-1}(p)| = f$.
\newline

Define the \textit{sort-number} of a permutation as the minimum number of times $SC_{231}$ must be applied to it in order for it to enter $\Av_n(132, 231)$.

Other than a computational analysis of $SC_{231}$ sort-numbers, our paper also improves the bounds for the maximum sort number in \cite{SiedelSun} to \[n-1 \leq m(n) \leq \frac{(n+1)(n-2)}{2}\] and proves additional properties of $SC_{231}$ mathematically.
\newline

A key idea for exploring this question is what we will call the \textit{index} of a permutation. Define the \textit{index} of a permutation as:
\begin{itemize}
\item $1$ if $1$ and $2$ are not adjacent.
\item $2$ if $1$ and $2$ are adjacent but $3$ is not adjacent to either one.
\item $3$ if $1$ and $2$ are adjacent, $3$ is adjacent to either one, but $4$ is not adjacent to $1, 2,$ or $3$. \\
\dots
\item $n$ if $1$ and $2$ are adjacent, $3$ is adjacent to either one, $4$ is adjacent to some entry in $1$ through $3$, and $5$ is adjacent to some entry in $1$ through $4$, ..., and $n$ is adjacent to some entry in $1$ through $n-1$.
\item Note that index $n-1$ does not exist, because if $1$ through $n-1$ are all in a block, then $n$, being the only remaining entry, must be adjacent to that block.
\end{itemize}
For example, $45231$ has index $1$, then $SC_{231}(45231) = 13254$ has index $1$, then $SC_{231}(13254) = 35421$ has index $2$, then $SC_{231}(35421) = 51243$ has index $2$, then $SC_{231}(51243) = 43215$ has index $5$.

A permutation is a periodic point if and only if its index is $n$, because both mean the permutation has no \textit{peaks}. We will define a \textit{peak} as an entry $i$ in a permutation $p$ that is not the first entry or the last entry, such that $p(i-1) < p(i) > p(i+1)$.

This idea is important because it is a monovariant for $SC_{231}$. When you apply $SC_{231}$, the permutation's index can stay the same or increase, but it can not decrease. This was proved as Claim 4.5 in \cite{SiedelSun}.
\newline

We will use Siedel and Sun's definitions of \textit{pre-popped} and \textit{post-popped}. While applying a stack such as $SC_{231}$, if an entry is popped from the stack when there are still entries in the input, then it is \textit{pre-popped}. Otherwise, let the entry be \textit{post-popped}.
\section{Methodology}

We used Java code to \textbf{compute} the sort-number of all permutations up to length $14$, further than Seidel and Sun's length $11$ counterexample. We also used Java code and the statistical t-test to \textbf{estimate} the average sort-number of length $n$ permutations for $15 \leq n \leq 1000$.
\newline\newline\noindent
Due to the length of the code (around 130 lines for each program), the full code will be provided upon request from the authors or be listed separate from the paper.
\newline\newline\noindent
The computation code is organized into sections:
\begin{itemize}
\item The main section contains a for loop so that the code processes all permutations between a start permutation and an end permutation (the user chooses these). It also uses a while loop to assist with determining the sort-number of the permutation.
\item The vshaped function takes a permutation $p$ and outputs a boolean representing whether the permutation is a periodic point of $SC_{231}$ (meaning no peaks).
\item The SC231 function takes in a permutation $p$ as an array and outputs $SC_{231}(p)$. This takes $O(n)$ runtime, where $n$ is the length of the permutation.
\item The nextPerm function takes in a permutation as an array and outputs the next permutation in this sequence: $1, 12, 21, 123, 132, 213, 231, 312, 321, 1234, \dots$.

This takes $\Theta(1)$ runtime because the number of operations is proportional to the length of the decreasing run at the end ($501432$ has a length $3$ decreasing run at the end, $1423$ has a length $1$ decreasing run at the end).
\end{itemize}

The output in the code shows the first permutation after the end permutation. The grid of numbers has $15$ rows and $30$ columns ($15$ and $30$ can be changed). The number in the $i$th row (top to bottom, starts at $1$) and $j$th column (left to right, starts at $0$) represents the number of permutations in the given range that have length $i$ and sort-number $j$.

The computation required approximately $75$ minutes to run for all permutations of length $\leq 13$ using $8\%$ of a $10$-core CPU. The code took around $5$ hours to run for all permutations of length $14$ using $56\%$ of a $10$-core CPU.
\newline\newline\noindent
The estimation code is organized into sections:
\begin{itemize}
\item The main section uses a for loop to loop through the desired permutation lengths, a for loop to generate the desired number of permutations for each length, and a while loop to compute each permutation's sort-number.
\item The vshaped function and the SC231 function are still the same.
\item The nextPerm function is replaced by the randPerm function, which uses the Fisher-Yates Shuffle Algorithm \cite{FisherYates} to generate a random permutation in $O(n)$ runtime.
\end{itemize}
Based on the structure of the code and the code's progress over time, we believe the runtime for estimating the average sort-number of permutations of length $n$ is $\Theta(n \cdot m \cdot a(n))$, where $m$ is the sample size and $a(n)$ is defined at the start of the Observations section.

The console output shows the progress of the program. It outputs the value $n$ it is currently processing. The first five rows of the file output represent the numbers up to $1000$, the mean of the sort-numbers for each length, the standard deviation of the sort-numbers for each length, the lower bounds of the 99\% confidence intervals, and the upper bounds of the 99\% confidence intervals, respectively. The remaining grid of numbers are all the calculated sort-numbers from the randomly generated permutations.

\section{Data Display}

Due to the length of the data (two $14 \times 30$ tables and four length $986$ lists), the full data will be provided upon request from the authors or be listed separate from the paper.

Tables 3 and 4 are from the computational code, while Table 2 and Graph 5 are from the estimation code.

Here is an overview of the data. For tables $3$ and $4$, the third column's full title is the number of permutations with the maximum sort-number. The average sort-number is computed as a Java double variable, which can reliably represent up to $15$ decimal digits.

\begin{figure}[H]
\centering\footnotesize\noindent
\begin{tabular}{l|l|l|l|l|l|l|l}
$n$ & 15 & 25 & 50 & 100 & 200 & 300 & 400 \\
\hline
Estimated Avg Sort-Number & 11.465 & 24.5025 & 69.3925 & 195.1575 & 541.89 & 964.975 & 1444.46 \\
Lower Bound of CI & 11.055417 & 23.870502 & 68.150269 & 192.79904 & 537.73894 & 959.12593 & 1436.5756 \\
Upper Bound of CI & 11.874583 & 25.134498 & 70.634731 & 197.51596 & 546.04106 & 970.82407 & 1452.3444 \\
\hline
$n$ & 500 & 600 & 700 & 800 & 900 & 1000 & \\
\hline
Estimated Avg Sort-Number & 1954.0025 & 2496.2875 & 3071.16 & 3676.2975 & 4291.2775 & 4920.8125 & \\
Lower Bound of CI & 1944.4278 & 2485.2875 & 3058.5745 & 3662.5979 & 4276.211 & 4904.9479 & \\
Upper Bound of CI & 1963.5772 & 2507.2875 & 3083.7455 & 3689.9971 & 4306.344 & 4936.6771 & \\
\end{tabular}
\caption*{\footnotesize Table 2: Estimated Average Sort-Numbers and Their 99\% Confidence Intervals}
\end{figure}

\begin{figure}[H]
\centering\small\noindent
\begin{tabular}{l|l|l|l}
$n$ & Max Sort-Number & $\#$ of Permutations & Average Sort-Number \\
\hline
1 & 0 & 1 & 0.0 \\
2 & 0 & 2 & 0.0 \\
3 & 2 & 1 & 0.5 \\
4 & 4 & 1 & 1.25 \\
5 & 6 & 2 & 2.1083333333333334 \\
6 & 8 & 1 & 2.948611111111111 \\
7 & 10 & 4 & 3.778373015873016 \\
8 & 12 & 2 & 4.629861111111111 \\
9 & 14 & 2 & 5.510821759259259 \\
10 & 16 & 88 & 6.427365244708994 \\
11 & 20 & 1 & 7.37919314674523 \\
12 & 21 & 351 & 8.366963456907033 \\
13 & 24 & 183 & 9.394762786403412 \\
14 & 28 & 2 & 10.465681418116624
\end{tabular}
\caption*{\footnotesize Table 3: Sort-Numbers of Permutations by Length}
\end{figure}

\begin{figure}[H]
\centering\small\noindent
\begin{tabular}{l|l|l|l}
Leading Entry & Max Sort-Number & $\#$ of Permutations & Average Sort-Number \\
\hline
1 & 26 & 3 & 11.93474462988786 \\
2 & 27 & 34 & 10.936350330161094 \\
3 & 24 & 130 & 10.111191946074758 \\
4 & 27 & 101 & 9.99864242865545 \\
5 & 27 & 93 & 10.093875898246559 \\
6 & 26 & 14 & 10.212733348666509 \\
7 & 28 & 2 & 10.307651192525325 \\
8 & 26 & 220 & 10.36966016718621 \\
9 & 27 & 12 & 10.408423886427359 \\
10 & 27 & 6 & 10.435903312383347 \\
11 & 27 & 21 & 10.452185017271823 \\
12 & 27 & 3 & 10.46205731109811 \\
13 & 27 & 11 & 10.469223845052838 \\
14 & 26 & 14 & 10.326896539995499
\end{tabular}
\caption*{\footnotesize Table 4: Sort-Numbers of Length $14$ Permutations by Leading Entry}
\end{figure}

\begin{figure}[H]
\centering
\includegraphics[width=0.45\linewidth]{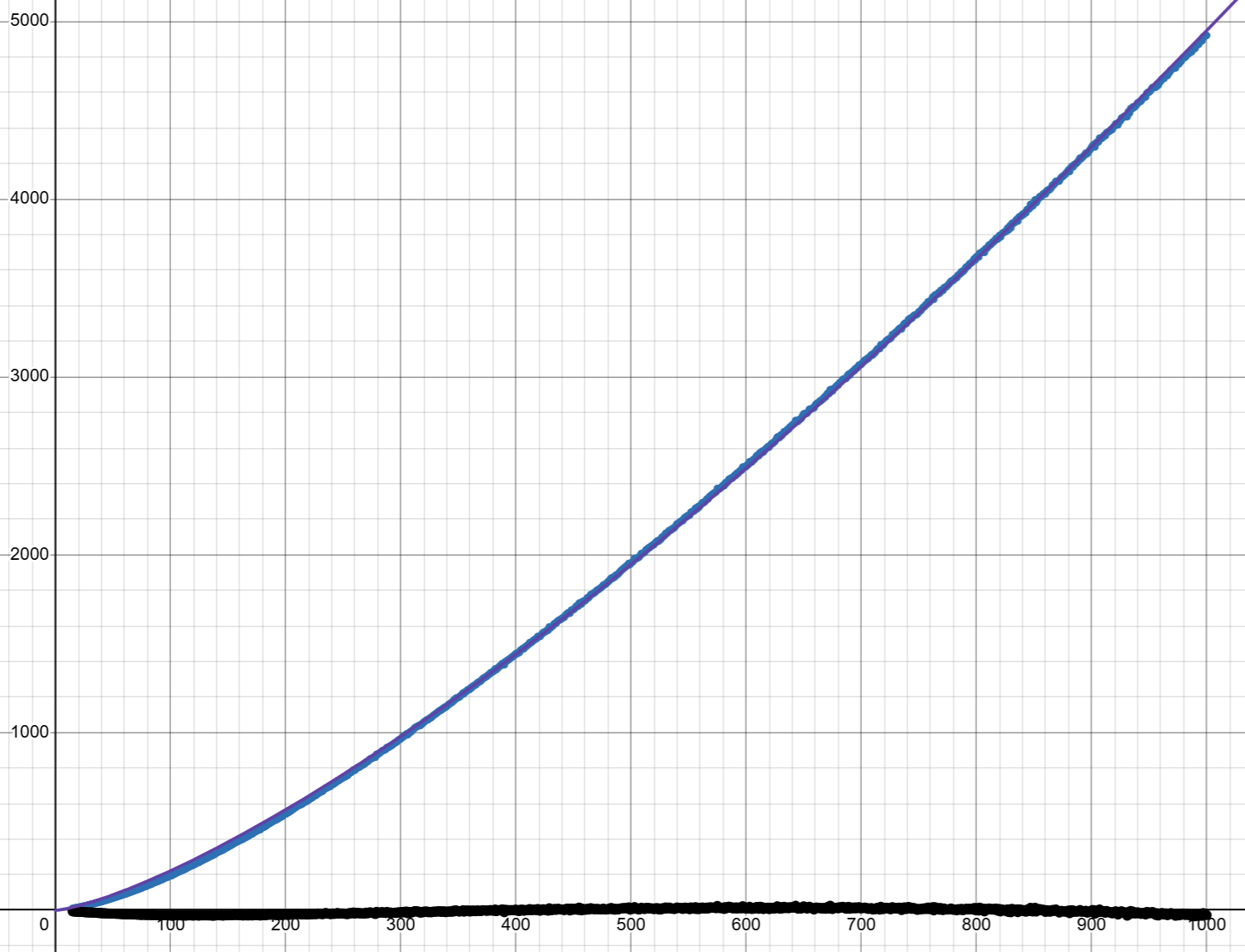}
\includegraphics[width=0.5\linewidth]{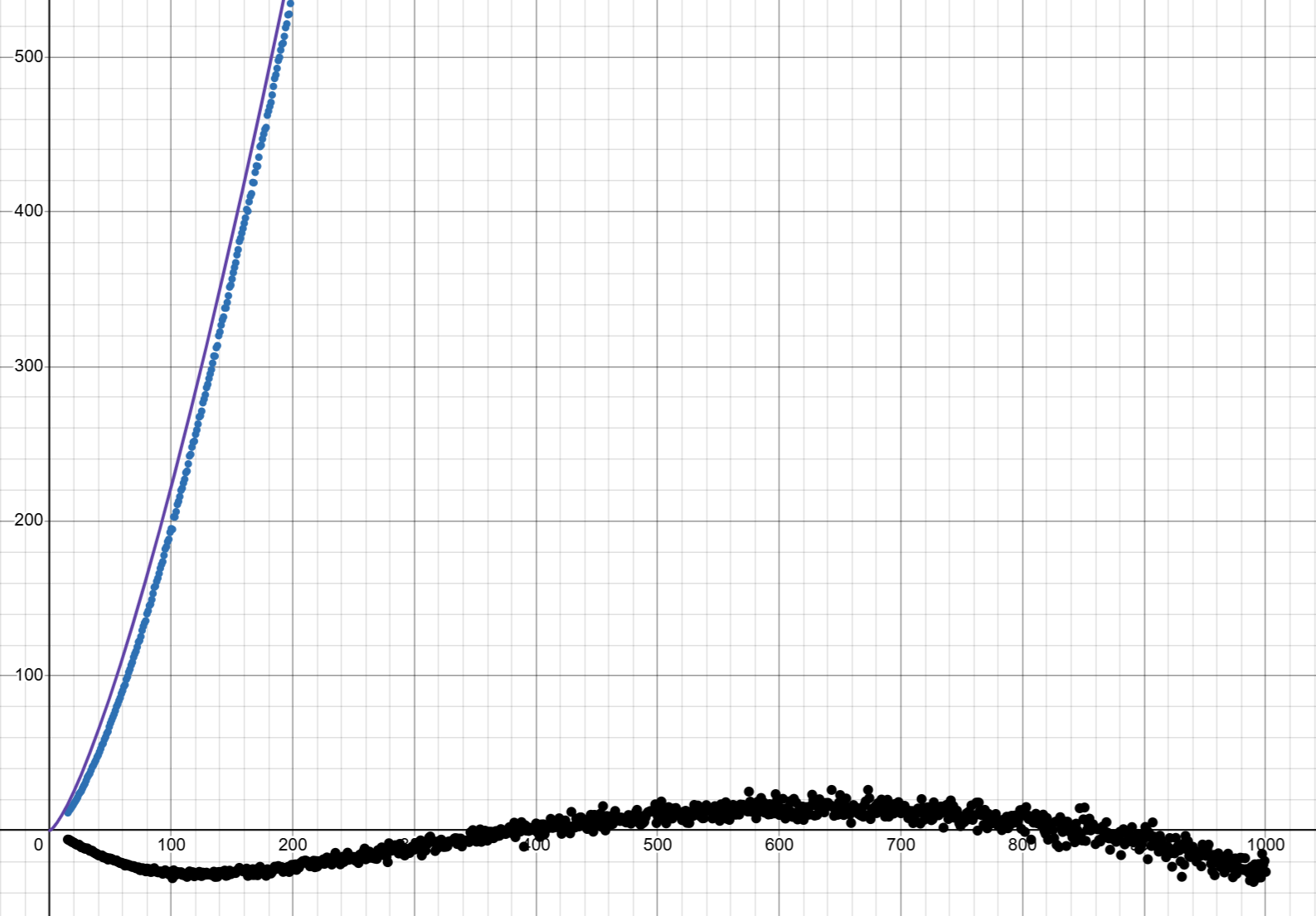}
\caption*{\footnotesize Graph 5: Average Sort-Number by Permutation Length with an $y = ax^b$ Curve Fit}
\end{figure}

\section{Observations}

Let $a(n)$ be the function that maps a positive integer $n$ to the average sort-number for permutations of length $n$.

Let $m(n)$ be the function that maps a positive integer $n$ to the maximum sort-number for permutations of length $n$.

Let $q(n, k)$ be the function that maps positive integers $n, k$ to the number of permutations of length $n$ with sort-number $k$.
\newline

The trend of $m(n)$ in Table $3$ suggests that the counterexample to Conjecture 1.1 in \cite{SiedelSun} --- the counterexample is a permutation of length $11$ and sort-number $19$ --- is not a lone counterexample but is likely the result of $m(n)$ having an above-linear trend.
\newline

The confidence intervals in Table $2$ are one sample t-intervals. Here are the conditions for using one sample t-intervals, all of which are satisfied here:
\begin{itemize}
\item Random Condition: Satisfied because we randomly chose the permutations among the permutations of length $n$.
\item Independence / 10\% Condition: Satisfied because we sampled with replacement.
\item Normal Condition: Satisfied. The sampling distribution for the average sort-number is approximately normal, because the population distribution (the sort-numbers of the permutations of length $n$) does not have strong skew or outliers based on the distributions for lengths $1$ through $14$, and the sample size is large ($400 \geq 30$).
\end{itemize}

When we graphed the permutation lengths with their average sort-numbers (calculated for $n \leq 14$ and estimated for $n \geq 15$) on Desmos \cite{Desmos}, we noticed that the residuals form a clear curve, meaning the $y=ax^b$ curve fit does not truly represent the data. However, it is still a good representation, because the correlation coefficient is $\geq 0.9999$. The values $a, b$ were $0.449796$ and $1.34712$.
\newline\newline\noindent
\textbf{Conjecture 4.1:} $\lim_{n\rightarrow\infty} \left(\frac{a(n)}{n}\right) > 0$. In other words, $a(n)$ is at least linear.
\begin{proof}
We outline partial progress below.

Taking inspiration from Lemma 4.3.1, one might think that because the permutations with index $1$ are pre-images of the permutations with index $> 1$, and the permutations with index $> 1$ are equivalent to the permutations of length $n-1$ (for example, $1243$ and $2143$ map to $132$. Each length $n-1$ permutation maps to exactly two length $n$ permutations with indices $> 1$, and no two length $n-1$ permutations share a length $n$ permutation), the average sort-number for the length $n$ permutations with index $1$ must be at least $1$ higher than the average sort-number for the permutations with length $n-1$.

This reasoning is mistaken because the length $n$ index $1$ permutations tend to map to the length $n-1$ permutations with below-average sort-numbers, as shown below for $n = 4$:

\begin{figure}[H]
\centering
\begin{asy}
label("\footnotesize1234", (0, 0), red);
label("\footnotesize2134", (0, -1), red);
label("\footnotesize3124", (0, -2), red);
label("\footnotesize3214", (0, -3), red);
label("\footnotesize4123", (0, -4), red);
label("\footnotesize4213", (0, -5), red);
label("\footnotesize4312", (0, -6), red);
label("\footnotesize4321", (0, -7), red);

draw((2.5,-1)--(0.5,-3), arrow=Arrow(TeXHead), deepgreen);
draw((2.5,-2.5)--(0.5,-5), arrow=Arrow(TeXHead), deepgreen);
draw((2.5,-3.5)--(0.5,-6), arrow=Arrow(TeXHead), brown);
draw((2.5,-5)--(0.5,-7), arrow=Arrow(TeXHead), brown);
draw((2.5,-6)--(0.5,-7), arrow=Arrow(TeXHead), deepgreen);
draw((2.5,-7)--(0.5,-7), arrow=Arrow(TeXHead), deepblue);

label("\footnotesize4132", (3, -1), green);
label("\footnotesize3142", (3, -2.5), green);
label("\footnotesize2143", (3, -3.5), red);
label("\footnotesize1243", (3, -5), red);
label("\footnotesize1423", (3, -6), green);
label("\footnotesize1432", (3, -7), blue);

draw((5.5,-0.5)--(3.5,-1), arrow=Arrow(TeXHead), deepgreen);
draw((5.5,-1.5)--(3.5,-1), arrow=Arrow(TeXHead), deepblue);
draw((5.5,-2.5)--(3.5,-2.5), arrow=Arrow(TeXHead), deepgreen);
draw((5.5,-3.5)--(3.5,-3.5), arrow=Arrow(TeXHead), brown);
draw((5.5,-5)--(3.5,-5), arrow=Arrow(TeXHead), brown);
draw((5.5,-6)--(3.5,-6), arrow=Arrow(TeXHead), deepgreen);
draw((5.5,-7)--(3.5,-7), arrow=Arrow(TeXHead), deepblue);

label("\footnotesize2314", (6, -0.5), green);
label("\footnotesize2431", (6, -1.5), blue);
label("\footnotesize2413", (6, -2.5), green);
label("\footnotesize3412", (6, -3.5), red);
label("\footnotesize3421", (6, -5), red);
label("\footnotesize3241", (6, -6), green);
label("\footnotesize2341", (6, -7), blue);

draw((8.5,-1.5)--(6.5,-1.5), arrow=Arrow(TeXHead), deepblue);
draw((8.5,-5)--(6.5,-5), arrow=Arrow(TeXHead), deepgreen);

label("\footnotesize1342", (9, -1.5), blue);
label("\footnotesize1324", (9, -5), green);

draw((11.5,-5)--(9.5,-5), arrow=Arrow(TeXHead), deepgreen);

label("\footnotesize4231", (12, -5), green);
\end{asy}
\caption*{\footnotesize Figure 6: Diagram of Length $4$ Permutations, Color-Coded by Distance Between Entries $1$ and $2$}
\end{figure}

The average sort-numbers of the index $>1$ permutations is $0.5$, but the average sort-numbers of the index $>1$ permutations with pre-images (weighted using the number of pre-images) is $0.4$.

Note that the diagram does not contradict Defant and Zheng's \cite{DefantZheng} Theorem 5.2, which states that for every $n \geq 2$, the max number of pre-images of a length $n$ permutation under $SC_{231}$ is $2^{n-2}$. Because $1234$ and $4321$ are both periodic points, the arrows $1234 \rightarrow 4321$ and $4321 \rightarrow 1234$ are not drawn.

Also, if we pair each length $n$ permutation with its complement (replacing each entry $i$ with $n+1-i$), the sum of the sort-numbers in the pair does not have to be linear in terms of $n$. For example, $15432$ has sort number $1$, while $51234$ has sort number $0$. Even if we extend the permutations to $198765432$ and $912345678$, the sum of the sort numbers is still $1$.

However, we did not use an approach that considers the number of pre-images.
\end{proof}

One may also extend Conjecture 4.1. From our data, $k$ is likely between $1$ and $1.3$.

\noindent
\textbf{Conjecture 4.1.2:} Find the value $k$ such that $\lim_{n\rightarrow\infty} \left(\frac{a(n)}{n^k}\right)$ is finite or prove that such a value does not exist.\newline

\noindent
\textbf{Theorem 4.2:} $m(n) \geq n-1$ for $n \geq 3$.

A lower bound for $m(n)$ was never provided in previous papers about $SC_{231}$, so we will provide one here.
\begin{proof}
Let $V_n$ be the permutation that is all the numbers with the same parity as $n$ listed in increasing order followed by all the numbers with the same parity as $n-1$ listed in decreasing order. For example, $V_3 = 132$ and $V_6 = 246531$.

We will prove that the sort-number of $\rev(V_n)$ is at least $n-1$ when $n \geq 3$. When the consecutive entries $n-1, n, n-2$ go into the stack, they form a \underline{132} pattern from top-to-bottom. Without a \underline{231} pattern in the stack, the stack simply reverses the permutation into $V_n$.

When we apply $SC_{231}$ to $V_n$, the entry $n$ gets popped from the stack when $n-1$ tries to enter the stack. The other entries are all left in the stack and are reversed when they exit the stack. The permutation is now $nV_{n-1}$.

We can do the same for $nV_{n-1}$. The entries $1$ and $2$ on the ends of $V_{n-1}$ cannot be pre-popped, since neither are large enough to be the $3$ in a \underline{231} pattern. The entries outside the $V_{n-1}$ can not re-enter between $1$ and $2$ since pre-popped entries are added to the front. Therefore, the entries outside $V_{n-1}$ do not affect how the elements inside $V_{n-1}$ are transformed by the $SC_{231}$ stack. We get that $nV_{n-1} = (n-1)V_{n-2}n$ because the only pre-popped entry inside $V_{n-1}$ is $n-1$.

This continues for $n-2, n-3, \dots, 3$. After $\dots V_3 \dots$ is sorted (using the $n-1$th $SC_{231}$), the entries $1$ and $2$ are at last adjacent. Before that, the presence of entries between $1$ and $2$ means there is a peak between $1$ and $2$.
\end{proof}

\noindent
\textbf{Theorem 4.3:} $m(n) \leq \frac{(n+1)(n-2)}{2}$ for $n \geq 3$.

Siedel and Sun \cite{SiedelSun} previously provided the upper bound $(n-1)(n-2)$. Through this claim, we will improve the upper bound.
\begin{proof}
First, we will provide a lemma that builds upon Claim 4.2:

\noindent
\textbf{Lemma 4.3.1:} In a permutation $p$, the number of entries between $1$ and $2$ decreases as $SC_{231}$ is repeatedly applied. The second, third, ... applications of $SC_{231}$ must decrease the number of entries between $1$ and $2$ by at least $1$ each time until $1$ and $2$ are next to each other.
\begin{subproof}
First of all, entries pre-popped from $SC_{231}$ are added to the front, not into the section between entries $1$ and $2$. Since $1$ and $2$ are never pre-popped (they are not large enough to be the $3$ in a $\underline{231}$ pattern), entries can only exit the section between $1$ and $2$. Therefore, the number of entries between $1$ and $2$ is decreasing (the first part of the lemma).

Let us say $p$ has entries between $1$ and $2$. Since $1$ and $2$ are both less than any entry in between them, there is a peak in between entries $1$ and $2$. The peak can either be a $\underline{132}$ pattern or a $\underline{231}$ pattern. If there exists a $\underline{132}$ pattern, applying $SC_{231}$ makes the $1$ enter the stack first, then the $3$; then trying to add the $2$ on top of the $3$ will cause the $3$ to be popped. After the middle entries of all $\underline{132}$ patterns are pre-popped, the remaining peaks must all be $\underline{231}$ patterns. The $SC_{231}$ reverses all the remaining entries, so these remaining peaks become $\underline{132}$ patterns. This shows that the second application of $SC_{231}$ must pre-pop at least one entry between $1$ and $2$ if there exists such an entry.

Applying this reasoning repeatedly proves the second part of the lemma.
\end{subproof}
From this lemma, we get a corollary:

\noindent
\textbf{Corollary 4.3.2:} For any permutation $p$ of index $1$, the permutation $SC_{231}^{n-1}(p)$ has index $>1$.
\begin{subproof}
There are at most $n-2$ entries between $1$ and $2$. From the lemma above, each application of $SC_{231}$ after the first removes at least one of these $n-2$ entries, so we need at most $n-1$ applications of $SC_{231}$ to remove all of these $n-2$ entries.
\end{subproof}
From Claim 4.5 of \cite{SiedelSun}, once the permutation has index $i$, we can treat it as a length $n+1-i$ permutation. Therefore, it takes at most $n-i$ applications of $SC_{231}$ to increase the index of an index $i$ permutation. Worst case scenario, we need $SC_{231}^{n-1}$ to increase the index of $p$ from $1$ to $2$, then $SC_{231}^{n-2}$ to increase the index from $2$ to $3$, ..., then $SC_{231}^2$ to increase the index from $n-2$ to $n$ (because index $n-1$ does not exist).
Therefore, \[m(n) \leq (n-1) + (n-2) + \dots + 2 = \frac{(n+1)(n-2)}{2}.\]
\end{proof}

\noindent
\textbf{Proposition 4.4:} $a(n)$ is at most quadratic.
\begin{proof}
Since $m(n) \geq a(n) \geq 0$, we can not have $a(n)$ rise faster than $m(n)$ in the long run. Therefore, this claim is proved directly using Claim 4.3, which gives us $m(n)$ is at most quadratic.
\end{proof}

\noindent
\textbf{Proposition 4.5:} $q(n, k)$ is not unimodal when $n$ is constant.
\begin{proof}
From our data, $q(4, 0) = 8$ and $q(4, 1) = 6$ and $q(4, 2) = 7$, breaking unimodality. This is not the only time unimodality is broken: for higher $n$, we also have $q(7, 2) = 1046$ and $q(7, 3) = 874$ and $q(7, 4) = 939$ among other examples.
\end{proof}

\noindent
\textbf{Proposition 4.6:} $q(n, k)$ is increasing when $k$ is constant. In fact, $q(n+1, k) \geq 2q(n, k)$.
\begin{proof}
We will map each permutation of length $n$ and sort-number $k$ to two unique permutations of length $n+1$ and sort-number $k$.

Using Claim 4.5 from \cite{SiedelSun}, the sort-number remains the same if we increase all the numbers in a length $n$ permutation by $1$, then replace $2$ with either $12$ or $21$. This map is reversible: take the $12$ or the $21$, contract it into a $2$, then decrease all the numbers by $1$.
\end{proof}

\noindent
\textbf{Claim 4.7:} When the leading entry is $1$ or $n$, the permutation's sort-number can not be exactly $2$.
\begin{proof}
Assume for the sake of contradiction there exists a permutation $p$ with leading entry $1$ and sort-number $2$. Since the leading entry can not be popped, $SC_{231}(p)$ has last entry $1$. For $p$ to have sort-number $2$ rather than $0$ or $1$, the entries in $SC_{231}(p)$ can not be in decreasing order; $SC_{231}(p)$ must have a peak.
\begin{itemize}
\item If there is a $\underline{231}$ pattern in $SC_{231}(p)$, since none of the three entries can be pre-popped (in order for an entry to be pre-popped, the two entries on either side must both be smaller), it will turn into a $\underline{132}$ pattern in $SC_{231}^2(p)$, a contradiction.
\item Otherwise, all peaks in $SC_{231}$ are $\underline{132}$ patterns, so they are all pre-popped and the remaining entries are reversed. But the $2$ and the $1$ in the $\underline{132}$ pattern are placed after the entry $1$ in $SC_{231}^2(p)$ in decreasing order, a contradiction.
\end{itemize}
Assume for the sake of contradiction there exists a permutation $p$ with leading entry $n$ and sort-number $2$. Since the leading entry can not be popped, $SC_{231}(p)$ has last entry $1$. For $p$ to have sort-number $2$ rather than $0$ or $1$, the entries in $SC_{231}(p)$ can not be in increasing order; $SC_{231}(p)$ must have a peak.
\begin{itemize}
\item If there is a $\underline{231}$ pattern in $SC_{231}(p)$, since none of the three entries can be pre-popped, it will turn into a $\underline{132}$ pattern in $SC_{231}^2(p)$, a contradiction.
\item Otherwise, all peaks in $SC_{231}(p)$ are $\underline{132}$ patterns, so they are all pre-popped and the remaining entries are reversed. But the pre-popped entries are placed before $n$ in $SC_{231}(p)$ and are smaller than $n$, a contradiction.
\end{itemize}
\end{proof}

\noindent
\textbf{Claim 4.8:} When the leading entry is $k \geq 1$ and $n \geq k$, there exists exactly $2^{k-2}$ (rounded to the nearest integer) permutations of length $n$ that are periodic points of $SC_{231}$.
\begin{proof}
When $k = 1$, the next entry must be higher, so the only way to avoid a peak is $123\dots n$.

When $k \geq 2$, the entry at the bottom must be $1$, and $k$ must be left of $1$. The entries $k+1, k+2, \dots, n$ must either be placed to the left of $k$ or be placed on the right. The former is impossible. The entries $2, 3, \dots, k-1$ must either be placed between $k$ and $1$ or to the right of $1$. Both are good, so there are $2^{k-2}$ ways to decide where these entries go. The order of the entries is automatic: the entries to the left of $1$ are in decreasing order, while the entries to the right of $1$ are in increasing order. So there are $2^{k-2}$ permutations of length $n$ that have no peaks.
\end{proof}

\noindent
\textbf{Conjecture 4.9:} $\lim_{n \rightarrow \infty} \left( \frac{m(n)}{n} \right) = \infty$. In other words, $m(n)$ is above linear.

While our computational data suggests $m(n)$ is above linear, we do not have a rigorous mathematical proof.

\begin{center}
\large
$\textbf{Acknowledgments}$
\end{center}

The author gratefully acknowledges Mikl\'os B\'ona for thought-provoking discussions and introducing the author to this research topic. The author also gratefully acknowledges Colin Defant for reviewing this paper and providing feedback.


\begin{thebibliography}{10}

\bibitem{Knuth}
Knuth, Donald E. The Art of Computer Programming, Volume 1: Fundamental Algorithms. Addison-Wesley, 1973.

\bibitem{West}
West, J. Permutations with restricted subsequences and stack-sortable permutations. Ph.D. Thesis, MIT, 1990.

\bibitem{CCF}
Cerbai, Giulio, et al. “Stack Sorting with Restricted Stacks.” \textit{Journal of Combinatorial Theory Series A}, vol. 173, July 2020, https://doi.org/10.1016/j.jcta.2020.105230.

\bibitem{DefantZheng}
Defant, Colin, and Kai Zheng. “Stack-Sorting with Consecutive-Pattern-Avoiding Stacks.” \textit{Advances in Applied Mathematics}, vol. 128, 18 Mar. 2021, https://doi.org/10.1016/j.aam.2021.102192. Accessed 24 June 2025.

\bibitem{SiedelSun}
Seidel, Ilaria, and Nathan Sun. “Periodic Points of Consecutive-Pattern-Avoiding Stack-Sorting Maps.” \textit{2308.05868}, arXiv, Aug. 2023, arxiv.org/pdf/2308.05868. Accessed 28 June 2025.

\bibitem{Zhao}
Zhao, William. “Stack-Sorting with Stacks Avoiding Vincular Patterns.” ArXiv.org, 2024, arxiv.org/abs/2410.17057. Accessed 30 Mar. 2026.

\bibitem{Kemeklis}
Kemeklis, Jurgis. “Fertility Numbers of Consecutive $S_3$ Pattern-Avoiding Stack-Sorting Maps.” ArXiv.org, 15 Sept. 2024, arxiv.org/abs/2408.05378. Accessed 2 Apr. 2026.

\bibitem{FisherYates}
GeeksforGeeks. “Shuffle a given Array Using Fisher–Yates Shuffle Algorithm.” GeeksforGeeks, 10 Oct. 2012, www.geeksforgeeks.org/dsa/shuffle-a-given-array-using-fisher-yates-shuffle-algorithm/.


\bibitem{Desmos}
Yi, Kai. “Sc231 Stats v2.” Desmos, 17 Apr. 2026, https://www.desmos.com/calculator/gamyae9s0y. Accessed 17 Apr. 2026.

\end{thebibliography}
\end{document}